\documentclass[12pt]{amsart}
\usepackage[margin=1in]{geometry}
\usepackage[english]{babel}
\usepackage[utf8]{inputenc}
\usepackage{subcaption}
\usepackage{amsmath}
\usepackage{amssymb}
\usepackage{amsfonts}
\usepackage{amsthm}
\usepackage{mathdots}
\usepackage{mathrsfs}
\usepackage[all]{xy}
\usepackage[pdftex]{graphicx}
\usepackage{color}
\usepackage{cite}
\usepackage{url}
\usepackage{indent first}
\usepackage[labelfont=bf,labelsep=period,justification=raggedright]{caption}
\usepackage[english]{babel}
\usepackage[utf8]{inputenc}
\usepackage{hyperref}
\usepackage[colorinlistoftodos]{todonotes}
\usepackage{tkz-fct}
\usepackage{tikz}
\usetikzlibrary{calc}
\usepackage{pgfplots}
\usepackage{multicol}
\PassOptionsToPackage{dvipsnames,svgnames}{xcolor}
\usepackage{textcomp}
\usepackage{bbm}
\usepackage[ruled,vlined]{algorithm2e}
\usepackage[shortlabels]{enumitem}

\newcommand{\QED}{\hspace{\stretch{1}} $\blacksquare$}

\theoremstyle{plain}
\newtheorem{thm}{Theorem}

\newtheorem{prop}[thm]{Proposition}

\theoremstyle{definition}

\theoremstyle{remark}

\numberwithin{equation}{section}
\numberwithin{thm}{section}

\allowdisplaybreaks

\begin{document}
	\title{An elementary proof of Ramanujan's identity for odd zeta values}
	\author{Sarth Chavan}
	\address{Euler Circle, Palo Alto, California 94306.}
	\email{sarth5002@outlook.com}
	\maketitle
	\begin{abstract}
	The main goal of this article is to present an elementary proof of Ramanujan's identity for odd zeta values. Our proof solely relies on a Mittag-Leffler type expansion for hyperbolic cotangent function and Euler's identity for even zeta values.
	\end{abstract}
\section{A Brief History and Introduction}
The Riemann zeta function $\zeta(s)$ is one of the most important special functions of Mathematics. While the critical strip $0 < \Re\left(s\right) < 1$ is undoubtedly the most important region in
the complex plane on account of the unsolved problem regarding location of non-trivial
zeros of $\zeta(s)$, namely, the Riemann Hypothesis, the right-half plane $\Re\left(s\right) > 1$ also has its
own share of interesting unsolved problems to contribute to.

It is quite well known that many number theoretic properties of odd zeta values are nowadays still unsolved mysteries, such as the rationality, transcendence and existence of closed-forms.
 Only in 1978 did Ap\'ery \cite{Apery} famously proved that $\zeta(3)$ is irrational. This was later reproved in a variety of ways by several authors, in particular Beukers \cite{Beukers} who devised a simple approach involving certain integrals over $[0, 1]^3$.
 
In the early 2000s, an
important work of Rivoal \cite{Rivol}, and Ball and Rivoal \cite{Brivol} determined that infinitely many values of $\zeta$ at odd
integers are irrational, and the work of Zudilin \cite{Zudilin} proved that at least one among $\zeta(5), \zeta(7), \zeta(9)$ and $\zeta(11)$
is irrational. 
A very recent result due to Rivoal and Zudilin \cite{RW1} states that at least
two of the numbers $\zeta(5), \zeta(7), \ldots , \zeta(69)$ are irrational.
Moreover, for any pair of positive integers $a$ and $b$, Haynes and Zudilin \cite[Theorem 1]{Haynes} have shown that either there are infinitely many $m \in \mathbb{N}$ for which $\zeta(am + b)$ is irrational, or the sequence $\{q_m\}_{m=1}^{\infty}$ of common denominators of the rational elements of the set $\{\zeta(a + b),\ldots, \zeta(am + b)\}$ grows super-exponentially, that is, $q_m^{1/m} \to \infty$ as $m \to \infty$.

Despite these advances, to this day no value of $\zeta(2n + 1)$ with $n\geqslant 2$ is known to be irrational. A folklore conjecture states that the numbers $\pi$,
$\zeta(3), \zeta(5), \zeta(7), \ldots$ are algebraically independent over the rationals. This conjecture is predicted by the Grothendieck’s period conjecture for mixed Tate motives. But both conjectures are far out of
reach and we do not even know the transcendence of a single odd zeta value. 

One should mention that Brown \cite{Brown} has in the past few years outlined a simple geometric approach to understand
the structures involved in Beukers’s proof of  irrationality of $\zeta(3) $ and how this may generalize to other odd zeta
values.

	Ramanujan made many beautiful and elegant discoveries in his short life of 32 years. One of the most remarkable formulas suggested by  Ramanujan  that has attracted the attention of several mathematicians over the years is the following intriguing identity involving the odd values of the Riemann zeta function \cite[1.2]{Berndt}:
\begin{thm}[Ramanujan's Formula for $ \zeta(2m+1) $]\label{MainThm}
	If $ \alpha $ and $ \beta $ are positive numbers such that $ \alpha\beta = \pi^2 $ and if $ m $ is a positive integer, then we have
	\begin{align}\alpha^{-m}\,\left\{\dfrac{1}{2}\,\zeta(2m + 1) + \sum_{n = 1}^{\infty}\dfrac{n^{-2m - 1}}{e^{2\alpha n} - 1}\right\} &-(- \beta)^{-m}\,\left\{\dfrac{1}{2}\,\zeta(2m + 1) + \sum_{n = 1}^{\infty}\dfrac{n^{-2m - 1}}{e^{2\beta n} - 1}\right\}\nonumber\\
	&=2^{2m}\sum_{k = 0}^{m + 1}\dfrac{\left(-1\right)^{k-1}B_{2k}\,B_{2m - 2k+2}}{\left(2k\right)!\left(2m -2k+2\right)!}\,\alpha^{m - k + 1}\beta^k. \label{maineq}\end{align}
	where $ B_m $ denotes the $ m $-th Bernoulli number with the exponential generating function $$\sum_{m\geqslant0}B_{m}\frac{z^m}{m!}=\dfrac{z}{\exp\left(z\right)-1}.$$
\end{thm}
Theorem \ref{MainThm} appears as Entry 21 in Chapter 14 of Ramanujan's second
notebook \cite[173]{Notebooks 2}.
The first published proof of Theorem \ref{MainThm} is due to S.L. Marulkar \cite{SL} although he was not aware that this
formula can be found in Ramanujan’s Notebooks. Grosswald too rediscovered this formula
and studied it more generally in \cite{Gross1}, \cite{Gross2}. Berndt \cite[Theorem 2.2]{BB1} derived a general formula
from which both Euler’s identity for even zeta values (1.1) and Ramanujan’s identity for odd zeta values  follow as special cases, thus showing that Euler’s and Ramanujan’s formulas are natural companions of each other.

Ramanujan's identity has a number of applications. For example, a contemporary interpretation of
the above identity, as given for instance in \cite{Encode}, is that it encodes fundamental transformation
properties of Eisenstein series on the full modular group and their Eichler integrals. This
observation is extended in \cite[Section 5]{BBAS} to weight $2k + 1$ Eisenstein series of level $2$ through
secant Dirichlet series. Moreover, Ramanujan’s identity also has applications in theoretical computer
science \cite{CS} in the analysis of special data structures and algorithms. More specifically, it
is used there to achieve certain distribution results on random variables related to 
data structures called \emph{tries}.
\vspace{0.03in}
\\
	B. Berndt and A. Straub in their article \cite[3.1]{Berndt} show that identity (\ref{maineq}) is equivalent to 
	\begin{align}\alpha^{-m}\sum_{n = 1}^{\infty}\dfrac{\coth\left(\alpha n\right)}{n^{2m+1}} &- (-\beta)^{-m}\sum_{n = 1}^{\infty}\dfrac{\coth\left(\beta n\right)}{n^{2m+1}}\nonumber \\&= - 2^{2m+1}\sum_{k = 0}^{m + 1}\dfrac{\left(-1\right)^{k}B_{2k}\,B_{2m + 2 - 2k}}{\left(2k\right)!\left(2m + 2 - 2k\right)!}\,\alpha^{m + 1 - k}\beta^k.\end{align}
Letting $ \alpha = \beta = \pi $ and replacing $ m $ with $ 2m + 1 $, we get the following \cite[3.2]{Berndt}:
	\begin{prop}
	\label{special_case}
	The following identity holds
	\begin{equation}\label{eqspecial}\sum_{n = 1}^{\infty}\dfrac{\coth\left(\pi n\right)}{n^{4m + 3}} = 2^{4m+2}\pi^{4m + 3}\sum_{k = 0}^{2m + 2}\dfrac{\left(-1\right)^{k+1}B_{2k}\,B_{4m + 4 -2k}}{\left(2k\right)!\left(4m+4-2k\right)!}\end{equation}
	where as before $ m \in \mathbb{N} $ and $ B_m $ denotes the $ m $-th Bernoulli number.
	\end{prop}
	This variation was apparently first established by M. Lerch \cite{Lerch}. Later proofs were given by G.N. Watson \cite{Watson}, H.F. Sandham \cite{Sadham}, J.R. Smart \cite{Smart} and F.P. Sayer \cite{Sayer}.
    
	The main purpose of this article is to present an elementary proof of Theorem \ref{MainThm} 
and Proposition \ref{special_case}. Our proof 
solely relies on a Mittag-Leffler type expansion for hyperbolic cotangent function and Euler's identity for even zeta values
\section{Preliminaries}
	\begin{prop}[Euler's Formula for $ \zeta(2m) $]
	\label{prop21}
		For $m \geqslant 1$, we have
		\begin{equation}\label{zeta2m}\zeta(2m) = \left(2\pi\right)^{2m}\dfrac{\left(-1\right)^{m - 1}B_{2m}}{2\left(2m\right)!}.\end{equation}
	\end{prop}
Euler's 
identity for even zeta values (\ref{zeta2m}) not only provides an elegant formula for evaluating $ \zeta(2m) $ for any integer $ m \geqslant 1 $, but it also gives us information about the arithmetical nature of even zeta values. For a fascinating account of the history and proof 
of this formula, we refer the reader to \cite{Apostol}. 
\begin{prop}[Mittag-Leffler expansion for hyperbolic cotangent function]
	\label{cor23}
	Let $ x \in \mathbb{C} \setminus{\{0\}} $, then the following identity holds
	\begin{equation}\label{ML2}
	\coth\left(\pi x\right)=\dfrac{1}{\pi x}+\dfrac{2x}{\pi}\sum_{k=1}^{\infty}\dfrac{1}{x^2+k^2}.
	\end{equation}
\end{prop}
Note that equation (\ref{ML2}) appears as identity 1.421.4 in \cite{GR}.
\section{Proof of Proposition \ref{special_case}}
	Replacing $ x $ with $ n $, dividing both sides of identity $ (\ref{ML2}) $ by $ n^{4m+2} $ and summing over $ n $ produces
	\[\dfrac{1}{2}\sum_{n = 1}^{\infty}\dfrac{\pi\coth(\pi n)}{n^{4m + 3}} - \dfrac{1}{2}\sum_{n = 1}^{\infty}\dfrac{1}{n^{4m + 4}} = \sum_{n = 1}^{\infty}\sum_{k = 1}^{\infty}\dfrac{1}{n^{4m + 2}\,(k^2 + n^2)}.\]
	Since the above double sum is absolutely convergent, we may switch the order of summation to get
	\[\dfrac{1}{2}\left(\sum_{n = 1}^{\infty}\dfrac{\pi\coth(\pi n)}{n^{4m + 3}} - \sum_{n = 1}^{\infty}\dfrac{1}{n^{4m + 4}}\right) = \sum_{k = 1}^{\infty}\sum_{n = 1}^{\infty}\dfrac{1}{n^{4m + 2}\,(k^2 + n^2)}.\]
	Next, we simply notice that
	\[
	\dfrac{1}{n^{4m + 2}\left(k^2 + n^2\right)} = 
	\dfrac{1}{n^{4m}k^2}\left(\dfrac{1}{n^2} - \dfrac{1}{k^2 + n^2}\right).\]
	Thus, we have
	\[ \sum_{n = 1}^{\infty}\sum_{k = 1}^{\infty}\dfrac{1}{n^{4m + 2}\left(k^2 + n^2\right)}= \zeta(2)\,\zeta(4m + 2) - \sum_{k = 1}^{\infty}\sum_{n = 1}^{\infty}\dfrac{1}{n^{4m}k^2\left(k^2 + n^2\right)}.\]
	Similarly, we find that
	\[\dfrac{1}{n^{4m}k^2\,(k^2 + n^2)}  =  
	\dfrac{n^2}{n^{4m}k^4}\left(\dfrac{1}{n^2} - \dfrac{1}{k^2 + n^2}\right).\] 
	Thus, we have
	\[\sum_{k = 1}^{\infty}\sum_{n = 1}^{\infty}\dfrac{1}{n^{4m}k^2\left(k^2 + n^2\right)}= \zeta(4)\,\zeta(4m) - \sum_{k = 1}^{\infty}\sum_{n = 1}^{\infty}\dfrac{1}{n^{4m - 2}k^4\left(k^2 + n^2\right)} .\]
	This recursive pattern that produces even zeta values continues as follows:
	\[\sum_{k = 1}^{\infty}\sum_{n = 1}^{\infty}\dfrac{1}{n^{4m - 2}k^4\left(k^2 + n^2\right)} = \zeta(6)\,\zeta(4m-2) - \sum_{k = 1}^{\infty}\sum_{n = 1}^{\infty}\dfrac{1}{n^{4m - 4}k^6\left(k^2 + n^2\right)}\]
	\[\sum_{k = 1}^{\infty}\sum_{n = 1}^{\infty}\dfrac{1}{n^{4m - 4}k^6\left(k^2 + n^2\right)} = \zeta(8)\,\zeta(4m-4) - \sum_{k = 1}^{\infty}\sum_{n = 1}^{\infty}\dfrac{1}{n^{4m - 6}k^8\left(k^2 + n^2\right)}\]
	and so on. At the end, we get
	\[\sum_{k = 1}^{\infty}\sum_{n = 1}^{\infty}\dfrac{1}{n^{2}k^{4m}\left(k^2 + n^2\right)} = \zeta(4m + 2)\,\zeta(2) - \sum_{k = 1}^{\infty}\sum_{n = 1}^{\infty}\dfrac{1}{k^{4m + 2}\left(k^2 + n^2\right)}.\]
	Therefore, we obtain the final result
	\begin{align}
	   \sum_{k = 1}^{\infty}\sum_{n = 1}^{\infty}\dfrac{1}{n^{4m + 2}\left(k^2 + n^2\right)}&= \zeta(2)\,\zeta(4m + 2) - \zeta(4)\,\zeta(4m) +  \cdots + \zeta(4m + 2)\,\zeta(2) \nonumber
	        \\
	    &- \sum_{k = 1}^{\infty}\sum_{n = 1}^{\infty}\dfrac{1}{k^{4m + 2}\left(k^2 + n^2\right)}\label{convoeq}
	\end{align}
	which can be easily proved by induction on $m$. Next, we observe that, by  symmetry,
	\[\sum_{k = 1}^{\infty}\sum_{n = 1}^{\infty}\dfrac{1}{n^{4m + 2}\left(k^2 + n^2\right)} = \sum_{k = 1}^{\infty}\sum_{n = 1}^{\infty}\dfrac{1}{k^{4m + 2}\left(k^2 + n^2\right)}.\]
	Therefore we have
	\begin{align*}2\sum_{k = 1}^{\infty}\sum_{n = 1}^{\infty}\dfrac{1}{n^{4m + 2}\left(k^2 + n^2\right)} &=  \left(-1\right)^m \zeta(2m+2)\,\zeta(2m+2) \\&+ 2\sum_{k=0}^{m-1} \left(-1\right)^k \zeta(4m+2-2k)\,\zeta(2k+2).\end{align*}
	Putting all things together produces
	\begin{align*}\sum_{n = 1}^{\infty}\dfrac{\pi\coth(\pi n)}{n^{4m + 3}}&= \zeta(4m + 4) + (-1)^m \zeta(2m+2)\,\zeta(2m+2) \\&+ 2\sum_{k=1}^m(-1)^{k-1}\,\zeta(4m-2k+4)\,\zeta(2k)\end{align*}
	where $k$ has been replaced with $k-1$ in the summand. Thus, it suffices to prove that
	\begin{align}
	\zeta(4m + 4) + \left(-1\right)^m &\zeta(2m+2)\,\zeta(2m+2) + 2\sum_{k=1}^m\left(-1\right)^{k-1}\zeta(4m-2k+4)\,\zeta(2k) \nonumber\\
	 = \,& \,2^{4m+2}\pi^{4m+4}\sum_{k = 0}^{2m+2}\dfrac{\left(-1\right)^{k-1}B_{2k}\,B_{4m+4-2k}}{\left(2k\right)!\left(4m+4-2k\right)!}.
	\end{align}
	After replacing  $m$ with $m - 1 $ in the above equation, it suffices to prove that
	\begin{align}
	\zeta(4m) + \left(-1\right)^{m-1}&\zeta(2m)\,\zeta(2m) + 2\sum_{k=1}^{m-1}\left(-1\right)^{k-1}\zeta(4m-2k)\,\zeta(2k)\nonumber
	  \\=\,&\,2^{4m-2}\pi^{4m}\sum_{k = 0}^{2m}\dfrac{\left(-1\right)^{k-1}B_{2k}\,B_{4m-2k}}{\left(2k\right)!\left(4m-2k\right)!}\label{Berno1}
	\end{align}
which is indeed just a matter of substitution; using Euler's  identity for $\zeta(2m)$, we have
\begin{align}
&\zeta(4m) + \left(-1\right)^{m-1}\zeta(2m)\,\zeta(2m) + 2\sum_{k=1}^{m-1}(-1)^{k-1}\,\zeta(4m-2k)\,\zeta(2k)
\nonumber\\&= 2^{4m-2}\pi^{4m}\left\{-\dfrac{B_{4m}}{\left(4m\right)!} + \dfrac{\left(-1\right)^{m-1}B_{2m}\,B_{2m}}{\left(2m\right)!\left(2m\right)!} +2\sum_{k=1}^{m-1}\dfrac{\left(-1\right)^{k-1}B_{2k}\,B_{4m-2k}}{\left(2k\right)!\left(4m-2k\right)!}\right\}\label{Berno2}
\end{align}
The summand in equation (\ref{Berno1}) is invariant under the substitution $ k \mapsto 2m - k $, so if we split it according to $k \in \{0, 2m\},\, k \in \{1,2,\ldots,m + 1, \ldots, 2m - 1\}$ and $k = m$, we get
\begin{equation}\sum_{k = 0}^{2m}\dfrac{\left(-1\right)^{k-1}B_{2k}\,B_{4m-2k}}{\left(2k\right)!\left(4m-2k\right)!}= -\dfrac{B_{4m}}{\left(4m\right)!} + \dfrac{\left(-1\right)^{m-1}B_{2m}\,B_{2m}}{\left(2m\right)!\left(2m\right)!} +2\sum_{k=1}^{m-1}\dfrac{\left(-1\right)^{k-1}B_{2k}\,B_{4m-2k}}{\left(2k\right)!\left(4m-2k\right)!}.\label{rfreq}\end{equation}
Combining equations (\ref{Berno1}), (\ref{Berno2}) and (\ref{rfreq}) gives us the desired result. \QED
\section{Proof of Theorem \ref{MainThm}}
Substituting $x\mapsto\left(\alpha k/\pi\right)$, with $\alpha\in\mathbb{R}_{>0}$, in identity (\ref{ML2}) yields
\begin{equation}\label{ML1}
	\dfrac{1}{e^{2\alpha k}-1}=-\dfrac{1}{2} + \dfrac{1}{2\alpha k} + \sum_{m=1}^{\infty}\dfrac{\alpha k}{\pi^{2}m^{2}+\alpha^{2}k^{2}}.
	\end{equation}
    Dividing both sides of the above identity by $\alpha^{n}k^{2n+1}$ and summing over $k$ produces
    \[\alpha^{-n}\left\{ \sum_{k=1}^{\infty}\dfrac{k^{-2n-1}}{e^{2\alpha k}-1}+\dfrac{1}{2}\,\zeta(2n+1)\right\} =\dfrac{\zeta(2n+2)}{2\alpha^{n+1}}+\sum_{k=1}^{\infty}\sum_{m=1}^{\infty}\dfrac{\alpha^{n+1}}{\left(\alpha^{2}k^{2}\right)^{n}\left(\pi^{2}m^{2}+\alpha^{2}k^{2}\right)}.\]  
 Next, we notice that, for an arbitrary positive
integer $n$ and two arbitrary non-vanishing sequences
$\{x_{k}\}$ and $\{z_{m}\}$, we have 
\begin{equation}\label{MM1}
\sum_{k=1}^{\infty}\sum_{m=1}^{\infty}\frac{1}{x_{k}^{n}\left(x_{k}+z_{m}\right)}=\sum_{k=1}^{\infty}\frac{1}{x_{k}^{n}}\sum_{m=1}^{\infty}\frac{1}{z_{m}}-\sum_{k=1}^{\infty}\sum_{m=1}^{\infty}\frac{1}{x_{k}^{n-1}z_{m}\left(x_{k}+z_{m}\right)}\end{equation}
assuming that the sums involved are convergent. Using the same recursive technique as in the previous section, we find that
\[\sum_{k=1}^{\infty}\sum_{m=1}^{\infty}\frac{1}{x_{k}^{n}\left(x_{k}+z_{m}\right)}=\sum_{k=1}^{\infty}\frac{1}{x_{k}^{n}}\sum_{m=1}^{\infty}\frac{1}{z_{m}}-\sum_{k=1}^{\infty}\frac{1}{x_{k}^{n-1}}\sum_{m=1}^{\infty}\frac{1}{z_{m}^{2}}\]
\begin{equation}\label{MM2}+\sum_{k=1}^{\infty}\frac{1}{x_{k}^{n-2}}\sum_{m=1}^{\infty}\frac{1}{z_{m}^{3}}+\cdots+(-1)^{n-1}\sum_{k=1}^{\infty}\frac{1}{x_{k}}\sum_{m=1}^{\infty}\frac{1}{z_{m}^{n}}+(-1)^{n}\sum_{k=1}^{\infty}\sum_{m=1}^{\infty}\frac{1}{z_{m}^{n}\left(x_{k}+z_{m}\right)}.
\end{equation}
For notational convenience, set 
\[
\sum_{k=1}^{\infty}\frac{1}{x_{k}^{n}}=\zeta_{x}(n),\quad\sum_{k=1}^{\infty}\frac{1}{z_{m}^{n}}=\zeta_{z}(n).
\]
These are also known as \emph{quasisymmetric zeta function} (the correct analog of a Riemann zeta function in the ring of quasisymmetric functions). 

Thus, equation  (\ref{MM2}) can now be rewritten as 
\begin{equation}\label{MM3}
\sum_{k=1}^{\infty}\sum_{m=1}^{\infty}\frac{1}{x_{k}^{n}\left(x_{k}+z_{m}\right)}=\sum_{p=0}^{n-1}(-1)^{p}\zeta_{x}(n-p)\,\zeta_{z}(p+1)+(-1)^{n}\sum_{k=1}^{\infty}\sum_{m=1}^{\infty}\frac{1}{z_{m}^{n}\left(x_{k}+z_{m}\right)}
\end{equation}
Substituting $x_{k}\mapsto\alpha^{2}k^{2}$ and $z_{m}\mapsto\pi^{2}m^{2}$
in identity (\ref{MM3}) produces
\begin{align}  \sum_{k=1}^{\infty}\sum_{m=1}^{\infty}\dfrac{1}{\left(\alpha^2k^2\right)^n\left(\alpha^2k^2 + \pi^2m^2\right)}
&=\sum_{p = 0}^{n-1}\left(-1\right)^{p}\dfrac{\zeta(2n-2p)}{\alpha^{2n-2p}}\,\dfrac{\zeta(2p + 2)}{\pi^{2p + 2}}\nonumber \\&+ \left(-1\right)^n\sum_{k=1}^{\infty}\sum_{m=1}^{\infty}\dfrac{1}{\left(\pi^2m^2\right)^n\left(\alpha^2k^2 + \pi^2m^2\right)}.\label{MM4}  
\end{align}
Euler's formula for $\zeta(2m)$ allows us to transform the zetas into Bernoullis as follows:
\[\sum_{p = 0}^{n-1}\left(-1\right)^{p}\dfrac{\zeta(2n-2p)}{\alpha^{2n-2p}}\dfrac{\zeta(2p + 2)}{\pi^{2p + 2}} = \left(-1\right)^n2^{2n}\sum_{p = 0}^{n-1}\left(\dfrac{\pi}{\alpha}\right)^{2n-2p}\dfrac{\left(-1\right)^{p-1}B_{2n - 2p}\,B_{2p + 2}}{\left(2n - 2p\right)!\left(2p + 2\right)!}.\]
Moreover, using $ \alpha\beta = \pi^2 $ and replacing $ p $ with $p - 1$ 
in the summand produces
\[\sum_{p = 0}^{n-1}\left(-1\right)^{p}\,\dfrac{\zeta(2n-2p)}{\alpha^{2n-2p}}\dfrac{\zeta(2p + 2)}{\pi^{2p + 2}} = \left(-1\right)^n2^{2n}\sum_{p = 1}^{n}\left(\dfrac{\beta}{\alpha}\right)^{n-p+1}\dfrac{\left(-1\right)^{p}B_{2n - 2p+2}\,B_{2p}}{\left(2n - 2p+2\right)!\left(2p\right)!}.\]
Substituting this expression in equation (\ref{MM4}) we get
\begin{align}
    \sum_{k=1}^{\infty}\sum_{m=1}^{\infty}\dfrac{1}{\left(\alpha^2k^2\right)^n\left(\alpha^2k^2 + \pi^2m^2\right)} & = \left(-1\right)^n2^{2n}\sum_{p = 1}^{n}\left(\dfrac{\beta}{\alpha}\right)^{n-p+1}\dfrac{\left(-1\right)^{p}B_{2n - 2p+2}\,B_{2p}}{\left(2n - 2p+2\right)!\left(2p\right)!}\nonumber
    \\&+ \left(-1\right)^n \sum_{k=1}^{\infty}\sum_{m=1}^{\infty}\dfrac{1}{\left(\pi^2m^2\right)^n\left(\alpha^2k^2 + \pi^2m^2\right)}\label{MM5}
\end{align}
 Putting all things together produces
 \[\alpha^{-n}\left\{\sum_{k = 1}^{\infty}\dfrac{k^{-2n-1}}{e^{2\alpha k}-1} + \dfrac{1}{2}\,\zeta(2n+1)\right\}=\left(-1\right)^n\alpha^{n+1}\sum_{k=1}^{\infty}\sum_{m=1}^{\infty}\dfrac{1}{(\pi^2m^2)^n(\alpha^2k^2 + \pi^2m^2)}\]
 \[+\, \dfrac{\zeta(2n+2)}{2\alpha^{n+1}} + \alpha^{n+1}\left\{= \left(-1\right)^n2^{2n}\sum_{p = 1}^{n}\left(\dfrac{\beta}{\alpha}\right)^{n-p+1}\dfrac{\left(-1\right)^{p}B_{2n - 2p+2}\,B_{2p}}{\left(2n - 2p+2\right)!\left(2p\right)!}\right\}.\]
 Next, we notice that the remaining zeta term
 \[\dfrac{\zeta(2n+2)}{2\alpha^{n+1}} = \dfrac{\left(-1\right)^n}{2\alpha^{n+1}}\dfrac{ \left(2\pi\right)^{2n+2}B_{2n+2}}{2\left(2n+2\right)!} = \dfrac{\left(-1\right)^n\beta^{n+1}2^{2n}B_{2n+2}}{\left(2n+2\right)!}\]
 is the missing $ p = 0 $-th term in the summand. 
 Thus, we have
 \begin{align}\alpha^{-n}\left\{ \sum_{k=1}^{\infty}\dfrac{k^{-2n-1}}{e^{2\alpha k}-1}+\dfrac{1}{2}\,\zeta(2n+1)\right\}  &- \alpha^{n+1}\left\{ \dfrac{\left(-1\right)^{n}}{\pi^{2n}}\sum_{k=1}^{\infty}\sum_{m=1}^{\infty}\dfrac{1}{m^{2n}\left(\alpha^{2}k^{2}+m^{2}\pi^{2}\right)}\right\}\nonumber
 \\&=\left(-1\right)^{n}2^{2n}\sum_{p=0}^{n}\dfrac{\left(-1\right)^{p}B_{2n-2p+2}\,B_{2p}}{\left(2n-2p+2\right)!\left(2p\right)!}\,\alpha^{p}\beta^{n-p+1}.\label{MM6}
 \end{align}
 Notice that, using $\alpha\beta=\pi^{2}$ we obtain
\[
\sum_{k=1}^{\infty}\sum_{m=1}^{\infty}\dfrac{1}{m^{2n}\left(\alpha^{2}k^{2}+m^{2}\pi^{2}\right)}=\dfrac{1}{\alpha}\sum_{k=1}^{\infty}\sum_{m=1}^{\infty}\dfrac{1}{m^{2n}\left(\alpha k^{2}+\beta m^{2}\right)}=\dfrac{1}{\alpha}\sum_{k=1}^{\infty}\sum_{m=1}^{\infty}\dfrac{1}{k^{2n}\left(\alpha^{2}m^{2}+\beta k^{2}\right)},
\]
where the last equality is obtained  by replacing $k$ by $m$
and $m$ by $k$ so that we get 
\begin{equation}\label{MM7}
\sum_{k=1}^{\infty}\sum_{m=1}^{\infty}\dfrac{1}{m^{2n}\left(\alpha^{2}k^{2}+m^{2}\pi^{2}\right)}=\dfrac{\beta}{\alpha}\sum_{k=1}^{\infty}\sum_{m=1}^{\infty}\dfrac{1}{k^{2n}\left(\beta^{2}k^{2}+m^{2}\pi^{2}\right)}
\end{equation}
Dividing both sides of identity (\ref{ML2}) by $k^{2n+1}$ and summing over $k$ produces
\begin{equation}\label{MM8}
\beta\sum_{k=1}^{\infty}\sum_{m=1}^{\infty}\dfrac{1}{k^{2n}\left(\beta^{2}k^{2}+m^{2}\pi^{2}\right)}=\sum_{k=1}^{\infty}\dfrac{k^{-2n-1}}{e^{2\beta k}-1}+\dfrac{1}{2}\,\zeta(2n+1)-\dfrac{1}{2\beta}\,\zeta(2n+2)
\end{equation}
where we have replaced $\alpha$ with $\beta$. Thus, we have
\[
\sum_{k=1}^{\infty}\sum_{m=1}^{\infty}\dfrac{1}{m^{2n}\left(\alpha^{2}k^{2}+m^{2}\pi^{2}\right)}=\dfrac{1}{\alpha}\left\{ \sum_{k=1}^{\infty}\dfrac{k^{-2n-1}}{e^{2\beta k}-1}+\dfrac{1}{2}\,\zeta(2n+1)-\dfrac{1}{2\beta}\,\zeta(2n+2)\right\}.
\]
Finally, substituting this expression in equation (\ref{MM6}) yields
\begin{align*}
 & \alpha^{-n}\left\{ \sum_{k=1}^{\infty}\dfrac{k^{-2n-1}}{e^{2\alpha k}-1}+\dfrac{1}{2}\,\zeta(2n+1)\right\} -\left(-1\right)^{n}2^{2n}\sum_{p=0}^{n}\dfrac{\left(-1\right)^{p}B_{2n-2p+2}\,B_{2p}}{\left(2n-2p+2\right)!\left(2p\right)!}\,\alpha^{p}\,\beta^{n-p+1}\\
 & =\dfrac{\left(-1\right)^{n}\alpha^{n}}{\pi^{2n}}\left\{ \sum_{k=1}^{\infty}\dfrac{k^{-2n-1}}{e^{2\beta k}-1}+\dfrac{1}{2}\,\zeta(2n+1)\right\} - \dfrac{\left(-1\right)^{n}\alpha^{n+1}}{\pi^{2n}}\left\{\dfrac{1}{2\alpha\beta}\,\zeta(2n+2)\right\}
\end{align*}
\[=(-\beta)^{-n}\left\{ \sum_{k=1}^{\infty}\dfrac{k^{-2n-1}}{e^{2\beta k}-1}+\dfrac{1}{2}\,\zeta(2n+1)\right\} -\dfrac{\left(-1\right)^{n}\alpha^{n+1}}{\pi^{2n+2}}\dfrac{\left(-1\right)^{n}2^{2n+2}\pi^{2n+2}B_{2n+2}}{4\left(2n+2\right)!}\]
Next, we notice that the remaining term: \vspace{0.05in}
\[
-\dfrac{\left(-1\right)^{n}\alpha^{n+1}}{\pi^{2n+2}}\dfrac{\left(-1\right)^{n}2^{2n+2}\pi^{2n+2}B_{2n+2}}{4\left(2n+2\right)!}=\dfrac{\left(-1\right)^{2n+1}2^{2n}B_{2n+2}}{\left(2n+2\right)!}\,\alpha^{n+1}
\]
is the $p=n+1$-st term in the summand. Therefore, we deduce that
\begin{align}
\alpha^{-n}\left\{ \sum_{k=1}^{\infty}\dfrac{k^{-2n-1}}{e^{2\alpha k}-1}+\dfrac{1}{2}\,\zeta(2n+1)\right\} &-\left(-\beta\right)^{-n}\left\{ \sum_{k=1}^{\infty}\dfrac{k^{-2n-1}}{e^{2\beta k}-1}+\dfrac{1}{2}\,\zeta(2n+1)\right\} 
\\
&
=\left(-1\right)^{n}2^{2n}\sum_{p=0}^{n+1}\dfrac{\left(-1\right)^{p}B_{2n-2p+2}\,B_{2p}}{\left(2n-2p+2\right)!\left(2p\right)!}\,\alpha^{p}\beta^{n-p+1}.\nonumber
\end{align}
Replacing $p$ with $n-p+1$ in the above identity gives us the desired result. \QED

\section{Acknowledgements}
The author would like to thank Christophe Vignat for his guidance and support throughout the completion of this work.
The author is very much thankful to the reviewer for his/her valuable comments towards the improvement of the manuscript. The author would also like to thank Bruce Berndt and Atul Dixit for endorsing him on arXiv.


\begin{thebibliography}{10}
\bibitem{Lost Notebook} G. E. Andrews and B.C. Berndt, Ramanujan\textquoteright s
Lost Notebook, Part IV, Springer, New York, 2013. 
\bibitem{Apery}  R. Ap\'ery, Irrationalit\'e de $\zeta(2)$ et $\zeta(3)$, Ast\'erisque \textbf{61} (1979), 11-13.
\bibitem{Apostol} T. M. Apostol. Another elementary proof of Euler’s formula for $\zeta(2n)$. Amer. Math. Monthly, 425-431.
\bibitem{Brivol} K. Ball and T. Rivoal, Irrationalit\'e d’une infinit\'e de valeurs de la fonction z$\hat{\textrm{e}}$ta aux entiers impairs (French), Invent. Math. \textbf{146} no. 1 (2001), 193–207.
\bibitem{Notebooks 2} B. C. Berndt, Ramanujan\textquoteright s
notebooks, Part II, Springer, New York, 1989. 
\bibitem{Berndt}B. C. Berndt and A. Straub, Ramanujan’s formula for $\zeta(2n + 1)$, Exploring the Riemann Zeta Function,
H. Montgomeryand, A. Nikeghbali and M. Rassias, Eds., Springer, Cham, 2017, 13-34.
\bibitem{Notebook 3}B. C. Berndt, Ramanujan\textquoteright s Notebooks,
Part III, Springer, New York, 1991. 
\bibitem{BB1}
B. C. Berndt, Modular transformations and generalizations of several formulae of Ramanujan, Rocky
Mountain J. Math. \textbf{7} (1977), 147-189.
\bibitem{BBAS}
B. C. Berndt and A. Straub, On a secant Dirichlet series and Eichler integrals of Eisenstein series, Math.
Z. 284 No. 3-4 (2016), 827–852.
\bibitem{Beukers} F. Beukers, A note on the irrationality of $\zeta(2)$ and $\zeta(3)$, Bull. London Math. Soc. 11:3 (1979), 268-272.
\bibitem{Brown} F. Brown, Irrationality proofs for zeta values, moduli spaces and dinner parties, Mosc. J. Comb. Number Theory, \textbf{6} No. 2-3 (2016), 102-165. 
\bibitem{seminar} S. Fischler, Irrationalit\'e de valeurs de $\hat{\textrm{z}}$eta (d’apr\'es Ap\'ery, Rivoal, $\ldots$) (S\'eminaire Bourbaki 2002-2003, expos\'e num\'ero 910, 17 November 2002) Ast\'erisque \textbf{294} (2004), 27-62.
\bibitem{GR}
I.S. Gradshteyn and I.M. Ryzhik, Table of Integrals, Series and Products, Edited by D. Zwillinger
and V.H. Moll, Eighth Edition, Academic Press, 2014.
\bibitem{Gross2}
E. Grosswald, Comments on some formulae of Ramanujan, Acta Arith. \textbf{21} (1972), 25-34.
\bibitem{Gross1}
E. Grosswald, Die Werte der Riemannschen Zetafunktion an ungeraden Argumentstellen,
Nachr. Akad. Wiss. G$\ddot{\textrm{o}}$ttinger Math.-Phys. Kl. II (1970), 9-13.
\bibitem{Encode}
 S. Gun, M. R. Murty and P. Rath, Transcendental values of certain Eichler integrals, Bull. London Math. Soc., \textbf{43} No. 5 (2011), 939-952. 
 \bibitem{Haynes}
 A. K. Haynes and W. Zudilin, Hankel determinants of zeta values, SIGMA Symmetry Integrability
Geom. Methods Appl. \textbf{11}, 101, (2015), 1-5.
 \bibitem{CS}
P. Kirschenhofer and H. Prodinger, On some applications of formulae of Ramanujan in the analysis of
algorithms, Mathematika \textbf{38} No. 1 (1991), 14-33.
 \bibitem{Lerch} M. Lerch, Sur la fonction $\zeta(s)$ pour valeurs impaires de l’argument, J. Sci. Math. Astron. pub. pelo Dr. F. Gomes Teixeira, Coimbra \textbf{14} (1901), 65-69.
\bibitem{SL} S.L. Malurkar, On the application of Herr Mellin’s integrals to some series, J. Indian Math. Soc. \textbf{16}
(1925/26), 130-138.
\bibitem{Rivol} T. Rivoal, La fonction $\hat{\textrm{z}}$eta de Riemann prend une infinit\'e de valeurs irrationnelles aux entiers impairs, C. R. Acad. Sci.
Paris S\'er. I Math. 331:4 (2000), 267-270.
\bibitem{RW1}
T. Rivoal, W. Zudilin A note on odd zeta values. S\'em. Lothar. Combin. \textbf{81} (2020), Art. B81b.
\bibitem{Sadham} H. F. Sandham. Some infinite series, Proc. Amer. Math. Soc. \textbf{5} (1954), 430–436.
\bibitem{Sayer} F. P. Sayer. The sums of certain series containing hyperbolic functions. Fibonacci Quart, \textbf{14}, 1976.
\bibitem{Smart} J.R. Smart, On the values of the Epstein zeta function, Glasgow Math. J. \textbf{14} (1973), 1-12.
\bibitem{Watson} G. N. Watson. Theorems stated by Ramanujan (ii): Theorems on summation of series. J. London Math. Soc. \textbf{3} (1928), 216-225.
\bibitem{Zudilin} W. Zudilin, One of the numbers $\zeta(5), \zeta(7), \zeta(9), \zeta(11)$ is irrational, Uspekhi Math. Nauk 56:4 (2001).
\end{thebibliography}
\end{document}